\documentclass[12pt]{amsart}

\newtheorem{thm}{Theorem}[section]
\newtheorem{prob}[thm]{Problem}

\theoremstyle{definition}

\theoremstyle{remark}

\renewcommand{\c}{{\mathfrak c}}

\newcommand{\cov}{\mathsf{cov}}
\newcommand{\add}{\mathsf{add}}
\newcommand{\cof}{\mathsf{cof}}

\newcommand{\non}{\mathsf{non}}

\renewcommand{\c}{\mathfrak{c}}

\newcommand{\R}{\mathbb{R}}

\newcommand{\Cite}[2]{{\cite[#1]{#2}}}

\newcommand{\bq}{\begin{quote}}
\newcommand{\eq}{\end{quote}}
\renewcommand{\O}{\mathcal{O}}
\newcommand{\B}{\mathcal{B}}
\newcommand{\BG}{\B_\Gamma}

\newcommand{\sone}{\mathsf{S}_1}

\newcommand{\NN}{{{}^{\naturals}\naturals}}
\newcommand{\naturals}{{\mathbb N}}
\newcommand{\N}{\naturals}

\newcommand{\by}[2]{\par\hfill\emph{#1}, #2}
\newcommand{\noemlby}[1]{\par\hfill\emph{#1}}

\newcommand{\CE}{\textsc{CE}}

\newcommand{\be}{\begin{enumerate}}
\newcommand{\ee}{\end{enumerate}}
\newcommand{\bi}{\begin{itemize}}
\newcommand{\ei}{\end{itemize}}
\renewcommand{\i}{\item}

\newcommand{\arx}[1]{\texttt{http://arxiv.org/abs/#1}}
\newcommand{\online}[1]{The paper is available online at
\bq
\texttt{#1}
\eq
}

\title[$\mathcal{SPM}$ Bulletin \textbf{6} (November 2003)]{%
$\mathcal{SPM}$ Bulletin\\[0.5cm]
Issue number 6: November 2003 \CE{}}

\begin{document}
\maketitle

\tableofcontents

\section{Editor's note}
After a long break, we are back with some very interesting research announcements and
an open problem which is one of the most difficult, long lasting, and important problems
in the field. 

A major change in this bulletin is that from now on it will usually not appear monthly, but
more close to quarterly. Special announcements (if urgent) will be made by text emails.
Contributions to the next issue are, as always, welcome.

{\small
\subsection*{Previous issues}
The first issues of this bulletin,
which contain general information (first issue),
basic definitions, research announcements, and open problems (all issues) are available online:
\be
\i First issue: \arx{math.GN/0301011}
\i Second issue: \arx{math.GN/0302062}
\i Third issue: \arx{math.GN/0303057}
\i Fourth issue: \arx{math.GN/0304087}
\i Fifth issue: \arx{math.GN/0305367}
\ee

\subsection*{Contributions}
Please submit your contributions (announcements, discussions, and open problems)
by e-mailing us. It is preferred to write them
in \LaTeX{}.
The authors are urged to use as standard notation as possible, or otherwise give
a reference to where the notation is explained.
Contributions to this bulletin would not require any transfer of copyright,
and material presented here can be published elsewhere.

\subsection*{Subscription}
To receive this bulletin (free) to your
e-mailbox, e-mail us.
}

\medskip

\by{Boaz Tsaban}{tsaban@math.huji.ac.il}

\hfill \texttt{http://www.cs.biu.ac.il/\~{}tsaban}

\section{Research announcements}

\subsection{The Topological Version of Fodor's Theorem}
The following purely topological generalization is given of
Fodor's theorem (also known as the ``pressing down lemma''):
Let $X$ be a locally compact, non-compact $T_2$ space such that
any two closed unbounded\footnote{A set is bounded if it has compact closure.}
(c\,u\,b) subsets of $X$ intersect; call $S
\subset X$ stationary if it meets every c\,u\,b in $X$. Then for
every neighbourhood assignment $U$ defined on a stationary
set $S$ there is a stationary subset $T \subset S$ such that
$$
\cap \{U(x) : x\in T\} \neq \emptyset.
$$

Just like the ``modern'' proof of Fodor's theorem, our proof
hinges on a notion of diagonal intersection of c\,u\,b's, definable
under some additional conditions.

We also use these results to present an (alas, only partial)
generalization to this framework of Solovay's celebrated
stationary set decomposition theorem.

\by{Istvan Juhasz}{juhasz@renyi.hu}
\noemlby{A.\ Szymanski}

\newcommand{\psmP}{\mbox{{\rm CPA$_{\rm prism}$}}}
\subsection{Continuous images of big sets and additivity of $s_0$ under \psmP}
We prove that the Covering Property Axiom \psmP,
which holds in the iterated perfect set model, implies the following 
facts. 
\begin{itemize}
\item 
There exists a family ${\mathcal G}$ of uniformly continuous functions from 
$\R$ to $[0,1]$ such that $|{\mathcal G}|=\omega_1$ and for every 
$S\in[\R]^\c$ there exists a $g\in{\mathcal G}$ with $g[S]=[0,1]$. 

\item
The additivity of the Marczewski ideal $s_0$ of is equal to $\omega_1<\c$.  
\end{itemize}

\by{Krzysztof Ciesielski}{K\_Cies@math.wvu.edu}
\by{Janusz Pawlikowski}{pawlikow@math.uni.wroc.pl}

\newcommand{\psmCGame}{{\rm CPA$_{\rm cube}^{\rm game}$}}
\subsection{Uncountable $\gamma$-sets under axiom CPA$_{\rm cube}^{\rm game}$}
In the paper we formulate a Covering Property Axiom \psmCGame,
which holds in the iterated perfect set model,
and show that it implies
the existence of uncountable strong $\gamma$-sets in $\R$
(which are strongly meager)
as well as uncountable $\gamma$-sets in $\R$ which 
are not strongly meager. These sets must be
of cardinality $\omega_1<\c$, since
every $\gamma$-set is universally null, while
\psmCGame\ implies that every universally null
has cardinality less than $\c=\omega_2$.

We will also show that \psmCGame\ implies the existence of
a partition of $\R$ into $\omega_1$ null compact sets.

\by{Krzysztof Ciesielski}{K\_Cies@math.wvu.edu}
\by{Andr\'es Mill\'an}{amillan@math.wvu.edu}
\by{Janusz Pawlikowski}{pawlikow@math.uni.wroc.pl}

\subsection{Ultrafilters with property (s)}
  A set X which is a subset of the Cantor set has property $(s)$ (Marczewski
(Spzilrajn)) iff for every perfect set $P$ there exists a perfect set $Q$ contained
in $P$ such that $Q$ is a subset of $X$ or $Q$ is disjoint from $X$. Suppose $U$ is a
nonprincipal ultrafilter on $\omega$. It is not difficult to see that if $U$ is
preserved by Sacks forcing, i.e., it generates an ultrafilter in the generic
extension after forcing with the partial order of perfect sets, then $U$ has
property $(s)$ in the ground model. It is known that selective ultrafilters or
even $P$-points are preserved by Sacks forcing. On the other hand (answering a
question raised by Hrusak) we show that assuming CH (or more generally MA for
countable posets) there exists an ultrafilter $U$ with property $(s)$ such that $U$ does
not generate an ultrafilter in any extension which adds a new subset of $\omega$.

\online{http://arXiv.org/abs/math/0310438}
\by{Arnold W.\ Miller}{miller@math.wisc.edu}

\subsection{There may be no Hausdorff ultrafilters}
An ultrafilter $U$ is Hausdorff if for any two functions $f,g$ mapping $\N$ to $\N$,
$f(U)=g(U)$ iff $f(n)=g(n)$ for $n$ in some $X$ in $U$. We will show that it is
consistent that there are no Hausdorff ultrafilters.

\online{http://arXiv.org/abs/math/0311064}

\by{Tomek Bartoszynski}{tomek@math.boisestate.edu}
\by{Saharon Shelah}{shelah@math.huji.ac.il}

\subsection{Proper forcing and rectangular Ramsey theorems}
I prove forcing preservation theorems for products of definable partial
orders preserving the cofinality of the meager or null ideal. Rectangular
Ramsey theorems for related ideals follow from the proofs.

\online{http://arXiv.org/abs/math/0311135}
\by{Jindrich Zapletal}{zapletal@math.ufl.edu}

\newcommand{\HN}{\mathcal{HN}}
\subsection{Cardinal characteristics of the ideal of Haar null sets}
We calculate the cardinal characteristics of the
$\sigma$-ideal $\HN(G)$ of Haar null subsets of a Polish
non-locally compact group $G$ with invariant metric and show that
$\cov(\HN(G))\le\mathfrak b\le\max\{\mathfrak
d,\non(\mathcal{N})\}\le\non(\HN(G))\le \cof(\HN(G))>\min\{\mathfrak
d,\non(\mathcal{N})\}$. If $G=\prod_{n\ge 0}G_n$ is the product of abelian
locally compact groups $G_n$, then $\add(\HN(G))=\add(\mathcal{N})$,
$\cov(\HN(G))=\min\{\mathfrak b,\cov(\mathcal{N})\}$,
$\non(\HN(G))=\max\{\mathfrak d,\non(\mathcal{N})\}$ and
$\cof(\HN(G))\ge\cof(\mathcal{N})$, where $\mathcal{N}$ is the ideal of Lebesgue
null subsets on the real line. Martin Axiom implies that
$\cof(\HN(G))>2^{\aleph_0}$ and hence $G$ contains a Haar null
subset of $G$ that cannot be enlarged to a Borel or projective
Haar null subset of $G$. This gives a negative (consistent) answer
to a question of S.Solecki. To obtain these estimates we show that
for a Polish non-locally compact group $G$ with invariant metric
the ideal $\HN(G)$ contains all $o$-bounded subsets (equivalently,
subsets with the small ball property) of $G$.

\by{Taras Banakh}{tbanakh@franko.lviv.ua}

\subsection{A note on generalized Egorov's theorem}
We prove that the following generalized version of Egorov's theorem
is independent from the ZFC axioms of the set theory.

Let $\{f_n\}_{n\in\omega}$, $f_n : \langle 0, 1\rangle\to R$,
be a sequence of functions (not necessarily measurable) converging
pointwise to zero for every $x\in \langle 0, 1\rangle $. Then for
every $\varepsilon >0$, there are a set $A\subseteq \langle 0, 1\rangle$
of Lebesgue outer measure $m^* > 1-\varepsilon$ and a sequence of
integers $\{ n_k\}_{k\in \omega}$ with $\{ f_{n_k}\}_{k\in \omega}$
converging uniformly on~$A$.     

\by{Tomasz Weiss}{tomaszweiss@go2.pl}

\section{Problem of the month}
The following is extracted from \cite{futurespm}.
Borel's Conjecture, which was proved to be consistent by Laver,
implies that each set of reals satisfying $\sone(\O,\O)$ (and the classes below it) is
countable.
From our point of view this means that there do not exist ZFC examples
of sets satisfying $\sone(\O,\O)$.
A set of reals $X$ is a \emph{$\sigma$-set} if each $G_\delta$ set in $X$
is also an $F_\sigma$ set in $X$.
In \cite{CBC} it is proved that every element of $\sone(\BG,\BG)$
is a $\sigma$-set.
According to a result of Miller \cite{Mil79Len},
it is consistent that every $\sigma$-set of
real numbers is countable. Thus, there do not exist uncountable
ZFC examples satisfying $\sone(\BG,\BG)$.
The situation for the other classes, though addressed by top 
experts, remains open. In particular, we have the following.

\begin{prob}[\cite{pawlikowskireclaw}, \Cite{Problem 45}{CBC}, \cite{ideals}]
Does there exist (in ZFC) an uncountable set of reals satisfying $\sone(\BG,\B)$?
\end{prob}
By \cite{CBC}, this is the same as asking whether it is consistent
that each uncountable set of reals can be mapped onto a dominating
subset of $\NN$ by a Borel function.
This is one of the major open problems in the field.


\begin{thebibliography}{00}
\bibitem{ideals}
T.\ Bartoszynski and B.\ Tsaban,
\emph{Hereditary topological diagonalizations and the Menger-Hurewicz Conjectures},
Proceedings of the American Mathematical Society,
to appear.
\arx{math.LO/0208224}

\bibitem{Mil79Len}
A.\ W.\ Miller,
\emph{On the length of Borel hierarchies},
Annals of Mathematical Logic \textbf{16} (1979),
233--267.

\bibitem{pawlikowskireclaw}
J.\ Pawlikowski and I.\ Rec{\l}aw,
\emph{Parametrized Cicho\'n's diagram and small sets},
Fundamenta Mathematicae \textbf{147} (1995),
135--155.

\bibitem{CBC}
M.\ Scheepers and B.\ Tsaban,
\emph{The combinatorics of Borel covers},
Topology and its Applications \textbf{121} (2002),
357--382.
\arx{math.GN/0302322}

\bibitem{futurespm}
B.\ Tsaban,
\emph{Selection principles in Mathematics: A milestone of open problems},
Note di Matematica,
to appear.

\end{thebibliography}
\end{document}